\documentclass[12pt]{article}
\usepackage{enumerate}
\usepackage{graphics}
\usepackage{amsmath}
\usepackage{array}
\usepackage{amssymb}
\newtheorem{thm}{Theorem}

 \numberwithin{equation}{section}
\begin{document}
\def\thefootnote{\fnsymbol{footnote}}
\begin{center}
\noindent \Large{ \textbf{On Generalized $q$-logistic
Distribution and Its Characterizations}\\}
\vskip.2cm\noindent \large{\textbf{Seema S.
Nair$^1$ and Nicy Sebastian$^2$\footnote{\textit{E-mail}: seema.cms@gmail.com, nicycms@gmail.com
}}}
 \vskip.05cm
 \noindent {\small{{\it $^1$ Department of Statistics,
 St. Gregorios College, Kottarakkara,
 India\\$^2$Department of Statistics,
 St. Thomas College, Thrissur,
 India }}}
\end{center}
 \vskip.2cm \noindent
{\bf\large{{Abstract}}}
Several generalizations of the logistic distribution, and certain
related models, are proposed by many authors for modeling various
random phenomena such as those encountered in data engineering,
pattern recognition, and reliability assessment studies. A
generalized $q$-logistic distribution is discussed here in the light
of pathway model, in which the new parameter $q$ allows increased
flexibility for modeling purpose. Also, we discuss different properties of the two generalizations of the $q-$logistic
distributions, which can be used to model the data exhibiting a unimodal density having
some skewness present. The first generalization is carried out using the basic idea of
Azzalini (1985) and we call it as the skew $q$-logistic distribution. It is observed that the
density function of the skew $q-$logistic distribution is always unimodal and log-concave
in nature.

\section{Introduction}
Logistic curves have been used as models in numerous applications. In
particular, because of its roles in analyzing bioassay and quantal
response experiments, a lot of research has been reported in
literature in studying the properties and applications of the
generalized logistic models. Mathai (2003) studied the distribution of order statistics from a logistic population, and pointed out some applications in survival analysis. In reliability theory, the classification of the tail behavior of
life distributions is helpful for characterizing their aging properties.
As it turns out, according to two classification systems
described in Rojo (1996), the logistic distribution is medium-tailed; that is, $\ln (1-F(x))$ is approximately linear in the tails, where $F(x)$
denotes the cumulative distribution function. The selection of an
appropriate model is paramount when studying the implications
for a system's integrity, or safety; or when attempting to identify
failure modes, and assessing future performance. The logistic
distribution allows generalization in many forms as could be seen in
George and Ojo (1980), Balakrishnan and Joshi (1983b), Balakrishnan and Leung (1988a), Wu et al.
(2000), Mathai (2003), Mathai and Provost (2006), Olapade (2006).

Because of the flexibility, much attention has been given to the
study of generalized models in recent times. The generalized model
proposed in this paper is referred to as $q$ analogue of the
logistic distribution in which the additional parameter $q$, called
as pathway parameter, is incorporated in its density function. In
addition to this, we have considered a generalized $q$-logistic
distribution by introducing location and scale parameters $\mu$ and
$\theta$ respectively which will
 result more flexible model than the standard $q$-logistic distribution.  The role of the
two additional parameters is to introduce skewness and to vary tail weights and provide greater
flexibility in the shape of the generalized distribution and consequently in modeling observed data.
It may be mentioned that although several skewed distribution functions exist on the positive real
axis, not many skewed distributions are available on the whole real line, which are easy to use
for data analysis purpose.
\section{Generalized $q$-logistic Model}
 Let $x$ be $q$-extended type-2 beta random variable having density
\begin{equation}\label{eqn2}
g(x)=c_1x^{\alpha-1}[1+a(q-1)x]^{-\frac{\beta}{q-1}},~a>0,~\beta>0,~x>0
\end{equation}
where $q>1$ is known as the pathway parameter through which one can
move from one functional form to another, see Mathai (2005), Mathai and Haubold (2007).

Suppose that we make the transformation  $y=\mu+\theta \ln x$, where
$x$ is distributed as in (\ref{eqn2}) for $q>1$, one has generalized
$q$-logistic density. That is
$$y=\mu+\theta \ln x\sim
G_qLD(\alpha,\beta,\mu,\theta,q),$$ and the corresponding density function has the following functional form
\begin{equation}\label{eqn1}
f(x)=\left\{ \begin{array}{ll} C [\rm{e}^{\frac{y-\mu}{\theta}}]
^\alpha[1+a(q-1)\rm{e}^{\frac{y-\mu}{\theta}}]^{-\frac{\beta}{q-1}},~q>1,-\infty<y<\infty,~\alpha>0,~a>0\\
0,~\rm{elsewhere},\end{array}\right.
\end{equation}
where
$$
C=\frac{[a(q-1)]^\alpha}{\theta}\frac{\Gamma(\frac{\beta}{q-1})}
{\Gamma(\alpha)\Gamma(\frac{\beta}{q-1}-\alpha)},~~\Re(\frac{\beta}{q-1}-\alpha)>0,~q>1.
$$ is the normalizing constant.
In particular when $q\rightarrow1$, one has
\begin{equation}
f(x)\rightarrow f_1(x)= C_1 [\rm{e}^{\frac{y-\mu}{\theta}}] ^\alpha
{\rm e}^{-a\beta\rm{e}^{\frac{y-\mu}{\theta}}},~~-\infty<y<\infty,~\alpha,~\beta,~\theta>0,~a>0,
\end{equation}
a generalized extreme value model, where $C_1=\frac{(a\beta)^\alpha}{\Gamma(\alpha)\theta}$ is the normalizing constant.

The density specified by (\ref{eqn1}) is plotted in the following
figures for $\alpha=2,~\beta=4,~a=1$ and different values of $q$.
The influence of additional parameters can be easily observed from
the following figures by putting different values to $\mu,\theta$
for each figure. It is observed that each type of generalized $q$-logistic curves established in this
 paper, the larger the pathway parameter $q$, the lower the mode of the corresponding distribution will be.
\begin{center}
\begin{eqnarray*}
\begin{array}{cc}
 \resizebox{6cm}{!}{\includegraphics{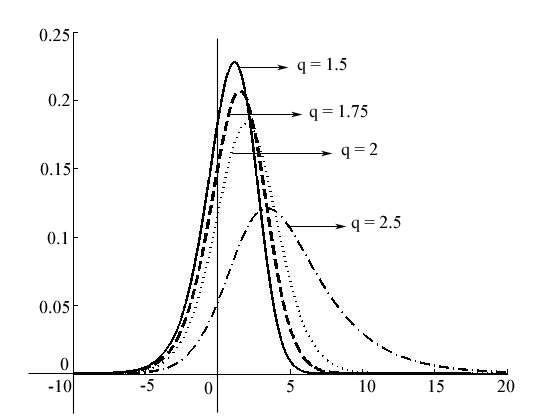}}~~~~&~~~~\resizebox{6cm}{!}{\includegraphics{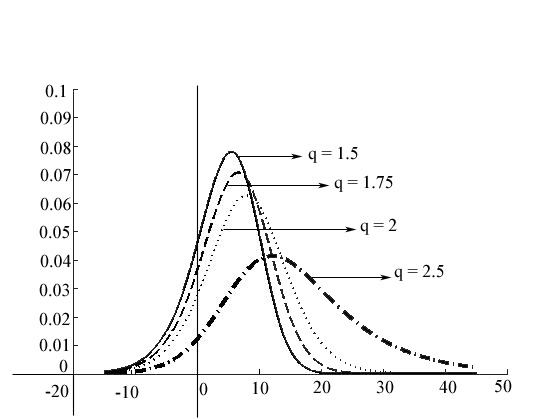}}\\
\text{Figure 1 {\scriptsize $GqL$ model for fixed $\mu=2$ and
$\theta=2.05$ }}~~~~~~~~~& \text{Figure 2 {\scriptsize $GqL$ model
for fixed $\mu=8$ and
$\theta=6$}}~~~
\end{array}
\end{eqnarray*}
\begin{eqnarray*}
\begin{array}{cc}
 \resizebox{6cm}{!}{\includegraphics{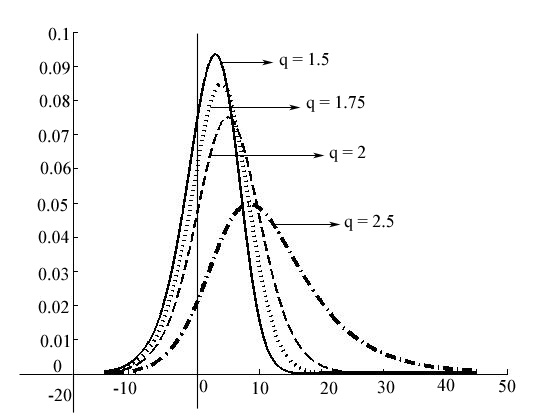}}~~~~&~~~ \resizebox{6cm}{!}{\includegraphics{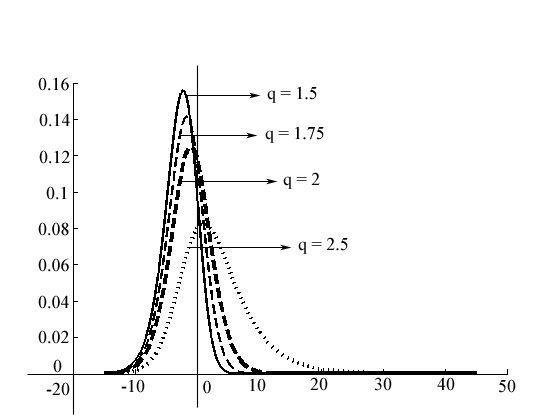}}\\
\text{Figure 3 {\scriptsize $GqL$ model for fixed $\mu=5$ and
$\theta=5$}}&
 \text{Figure 4 {\scriptsize $GqL$ model for
fixed $\mu=-1$ and $\theta=3$}}
\end{array}
\end{eqnarray*}
\end{center}

\section{Distribution and Survival Functions}
\noindent The distribution function of the model in (\ref{eqn1}) is
\begin{eqnarray}
F(t)&&=\int_{-\infty}^t f(y){\rm d}y\nonumber\\
&&=\frac{[a(q-1)]^\alpha}{\theta}\frac{\Gamma(\frac{\beta}{q-1})}
{\Gamma(\alpha)\Gamma(\frac{\beta}{q-1}-\alpha)}\int_{-\infty}^t
[\rm{e}^{\frac{y-\mu}{\theta}}]
^\alpha[1+a(q-1)\rm{e}^{\frac{y-\mu}{\theta}}]^{-\frac{\beta}{q-1}}
{\rm d}y,~~q>1\nonumber\\
&&=\frac{\Gamma(\frac{\beta}{q-1})}
{\Gamma(\alpha)\Gamma(\frac{\beta}{q-1}-\alpha)}\frac{[a(q-1)]^\alpha[{\rm
e}^{\frac{y-\mu}{\theta}}]^\alpha}{\alpha}\nonumber\\
&&~~~~~~~~~~~~\times{_2F_1}[\alpha,\frac{\beta}{q-1};(\alpha+1);-a(q-1){\rm
e}^{\frac{y-\mu}{\theta}}],~\text{when}~~q>1,~|a(q-1){\rm
e}^{\frac{y-\mu}{\theta}}|\leq 1.\nonumber
\end{eqnarray}
\noindent As a special case, when $\alpha=1,~q>1,~a>0,$ one has
\begin{eqnarray}
F_1(y)&=&\frac{\Gamma(\frac{\beta}{q-1})}
{\Gamma(\frac{\beta}{q-1}-1)}\int_{0}^{a(q-1){\rm
e}^{\frac{y-\mu}{\theta}}} [1+u]^{-\frac{\beta}{q-1}}
{\rm d}u,~~q>1,~0<a(q-1){\rm
e}^{\frac{y-\mu}{\theta}}<1\nonumber\\
&&=1-[1+a(q-1){\rm
e}^{\frac{y-\mu}{\theta}}]^{-(\frac{\beta}{q-1}-1)},~~|a(q-1){\rm
e}^{\frac{y-\mu}{\theta}}|<1
\end{eqnarray}
Hence in this case, the survival function is,
\begin{equation}
 \bar{F}_1(y)=[1+a(q-1){\rm
e}^{\frac{y-\mu}{\theta}}]^{-(\frac{\beta}{q-1}-1)},~~|a(q-1){\rm
e}^{\frac{y-\mu}{\theta}}|<1,
\end{equation}
which is used to model the life lengths of certain components of
interest in a device or system. The instantaneous failure rate or
the hazard rate, $\mu(t)$ of (\ref{eqn1}) can be given as
\begin{equation}
\mu(t)=\frac{f(t)}{1-F(t)}=\frac{{\rm d}}{{\rm d}t}\{-\ln(1-F(t))\}=c~\rm{e}^{\frac{t-\mu}{\theta}}[1+a(q-1){\rm
e}^{\frac{t-\mu}{\theta}}]^{-1},
\end{equation}
which can be shown to be an increasing function of $t$ and therefore
the model implies the aging effect. Also the model in (\ref{eqn1})
belongs to the family of IFR (Increasing Failure Rate) distributions
since $-\ln\{[1+a(q-1){\rm
e}^{\frac{y-\mu}{\theta}}]^{-(\frac{\beta}{q-1}-1)}\}$ is convex in
$t$, which is shown in figure 6.
\begin{center}
\begin{eqnarray*}
\begin{array}{cc}
 \resizebox{6cm}{!}{\includegraphics{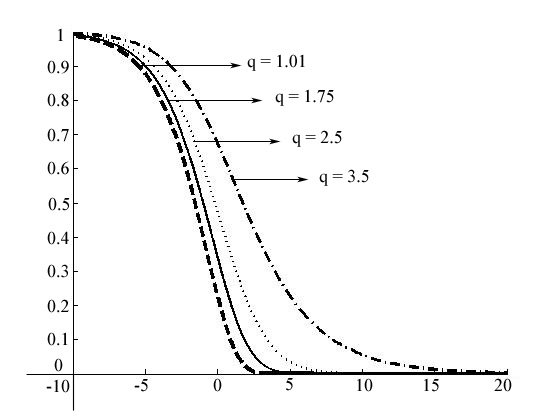}}~~~~~~~&~~~~~~~ \resizebox{6cm}{!}{\includegraphics{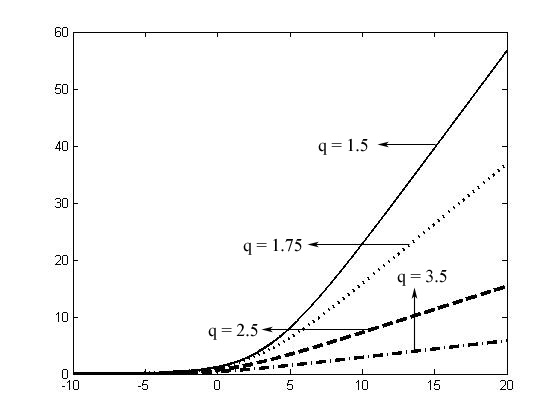}}\\
 \text{Figure 5: {\scriptsize The survival function of GqL }}&
 \text{Figure 6: {\scriptsize $-\ln(1-F_1(x))$}}\\
 \end{array}
\end{eqnarray*}
\end{center}
\begin{center}
 \resizebox{8cm}{5cm}{\includegraphics{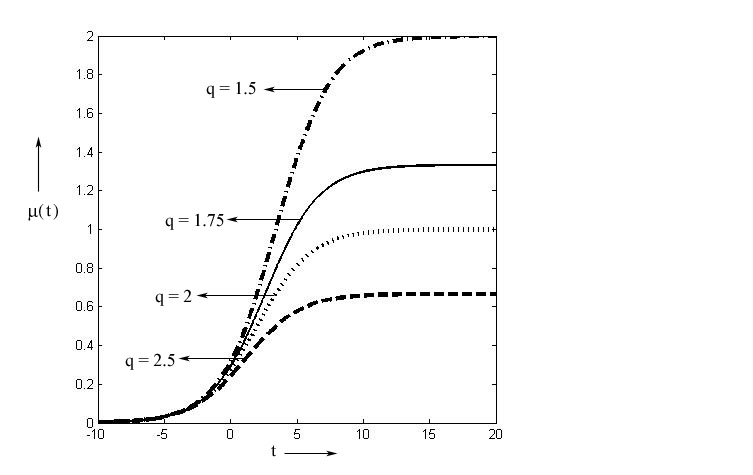}}
 \text{Figure 7: {\scriptsize The instantaneous failure rate $\mu(t)$ of GqL for various $q$. }}
\end{center}
\section{Characterization theorems based on generalized $q$-logistic distribution}
In this section, some theorems that characterize the generalized
$q$-logistic distribution are stated and proved.
\begin{thm}
Let $x$ be a continuously distributed random variable with density
function $f(x)$. Then the random variable $y=\mu+\theta
\ln[\frac{{\rm e}^x-1}{a(q-1)}]$ is a generalized $q$-logistic random
variable with parameters $(1,q)$ if and only if $x$ follows an
exponential distribution with parameter $\frac{\beta}{q-1},~q>1$.
\end{thm}
{\bf{Proof:}}\hskip.5cm Suppose that $x$ has exponential
distribution with parameter $\frac{\beta}{q-1}$ with density
function as
\begin{equation}\label{eqn5}
f(x)=\frac{\beta}{q-1}{\rm
e}^{-(\frac{\beta}{q-1})x},~~x>0,~\beta>0,~q>1
\end{equation}
Now let us take the transformation $$y=\mu+\theta \ln\left[\frac{{\rm
e}^x-1}{a(q-1)}\right]\Rightarrow\frac{y-\mu}{\theta}=\ln\left[\frac{{\rm
e}^x-1}{a(q-1)}\right]$$
$$\Rightarrow x=\ln[1+a(q-1)\rm{e}^{\frac{y-\mu}{\theta}}]\Rightarrow \frac{{\rm d}x}
{{\rm
d}y}=\frac{a(q-1)}{\theta}\rm{e}^{\frac{y-\mu}{\theta}}[1+a(q-1)\rm{e}^{\frac{y-\mu}{\theta}}]^{-1}$$
Then by Jaccobian of transformation, we find that
\begin{equation}
f(y)=\frac{a\beta}{\theta}\rm{e}^{\frac{y-\mu}{\theta}}
[1+a(q-1)\rm{e}^{\frac{y-\mu}{\theta}}]^{-(\frac{\beta}{q-1}+1)},~~-\infty<
y<\infty
\end{equation}
which is the probability density function of a generalized
$q$-logistic random variable $y$ with parameters
$(\frac{\beta}{q-1},1)$.

Conversely, suppose that $y$ is a generalized $q$-logistic random
variable, then the characteristic function of $y$ is given as
\begin{eqnarray}
\Phi_y(t)&=&a\beta\int_0^\infty {\rm e}^{it(\mu+\theta
\text{ln}u)}[1+a(q-1)u]^{-(\frac{\beta}{q-1}+1)}{\rm d}u,~~q>1\nonumber\\
&=&\frac{a\beta{\rm e}^{it\mu}}{[a(q-1)]^{it\theta+1}}\int_0^\infty z^{it\theta}
[1+z]^{-(\frac{\beta}{q-1}+1)}{\rm d}z,~~z=a(q-1)u\nonumber\\
&=&\frac{{\rm
e}^{it\mu}}{[a(q-1)]^{it\theta}}\frac{\Gamma(1+it\theta)
\Gamma(\frac{\beta}{q-1}-it\theta)}{\Gamma(\frac{\beta}{q-1})},~~\Re(\frac{\beta}{q-1}-it\theta)>0
\end{eqnarray}
$\text{for}~~q>1,~\Re(
\frac{\beta}{q-1}-it\theta)>0$. And also it is found that
\begin{equation}\label{eqn4}
\frac{{\rm e}^{it\mu}}{[a(q-1)]^{it\theta}}\frac{\Gamma(1+it\theta)
\Gamma(\frac{\beta}{q-1}-it\theta)}{\Gamma(\frac{\beta}{q-1})}={\rm
e}^{it\mu}\int_0^\infty\bigg[\frac{{\rm
e}^x-1}{a(q-1)}\bigg]^{it\theta}f(x){\rm d}x
\end{equation}
The only density function $f(x)$ satisfying the equation
(\ref{eqn4}) is the exponential distribution given in equation
(\ref{eqn5}). Hence the theorem is proved.

\begin{thm}
Suppose a continuously distributed random variable $x$ has a
$t$-distribution with $m$ degrees of freedom. Then the random
variable $y=\mu+\theta \ln[\frac{x^2}{a m (q-1)}]$ is
distributed according to generalized $q$-logistic random variable
with parameters $(\frac{1}{2},\frac{m+1}{2})$.
\end{thm}
{\bf{Proof:}}\hskip.5cm A random variable $x$ has a $t$-distribution
with $m$ degrees of freedom if
\begin{equation}
f(x)=\frac{\Gamma(\frac{m+1}{2})}{\Gamma(\frac{m}{2})\sqrt{(\pi
m)}}(1+\frac{x^2}{m})^{-(\frac{m+1}{2})},~~-\infty<x<\infty
\end{equation}
In order to apply one to one transformation, we split the range of variation of $x$ in to two, so that the density can be written as
\begin{equation}\label{eqn6}
f(x)=\left\{ \begin{array}{ll} f_1(x),~~-\infty<x<0\\
f_2(x),~~0<x<\infty\end{array}\right.
\end{equation}
where
\begin{equation}\label{eqn7}
f_1(x)=\left\{ \begin{array}{ll}\frac{\Gamma(\frac{m+1}{2})}{\Gamma(\frac{m}{2})\sqrt{(\pi
m)}}(1+\frac{x^2}{m})^{-(\frac{m+1}{2})},~~-\infty<x<0\\
0,~~otherwise\end{array}\right.
\end{equation}
and
\begin{equation}\label{eqn8}
f_2(x)=\left\{ \begin{array}{ll}\frac{\Gamma(\frac{m+1}{2})}{\Gamma(\frac{m}{2})\sqrt{(\pi
m)}}(1+\frac{x^2}{m})^{-(\frac{m+1}{2})},~~0<x<\infty\\
0,~~otherwise\end{array}\right.
\end{equation}
For $0<x<\infty$, suppose we take the transformation $y=\mu+\theta
\text{ln}[\frac{x^2}{a m (q-1)}]$, then $x=\sqrt{am(q-1)}{\rm
e}^{\frac{y-\mu}{2\theta}}$. Therefore, $\frac{{\rm
d}x}{{\rm dy}}=\frac{\sqrt{am(q-1)}}{2\theta}{\rm
e}^{\frac{y-\mu}{2\theta}}.$ Hence,
\begin{eqnarray}
f_1(y)&=&f_1(x)|_{y} \bigg|\frac{{\rm
d}x}{{\rm dy}}\bigg|\nonumber\\
&=&\frac{\Gamma(\frac{m+1}{2})}{\Gamma(\frac{m}{2})\Gamma(\frac{1}{2})}\frac{\sqrt{a(q-1)}}{2\theta}
{\rm e}^{\frac{y-\mu}{2\theta}} [1+a(q-1){\rm
e}^{\frac{y-\mu}{\theta}}]^{-\frac{m+1}{2}}
\end{eqnarray}
Since $f(x)$ is symmetric about zero, so it is clear that the transformed function $(f_2(y))$ is same for $-\infty<x<0$. Hence
\begin{eqnarray}
f(y)&=&\sum_{i=1}^2f_i(y)\nonumber\\
&=&\frac{\Gamma(\frac{m+1}{2})}{\Gamma(\frac{m}{2})\Gamma(\frac{1}{2})}\frac{\sqrt{a(q-1)}}{\theta}
{\rm e}^{\frac{y-\mu}{2\theta}} [1+a(q-1){\rm
e}^{\frac{y-\mu}{\theta}}]^{-\frac{m+1}{2}}
\end{eqnarray}

which is the probability density function for generalized
$q$-logistic random variables $(\frac{1}{2},\frac{m}{2})$.
Conversely, if $y$ is a generalized $q$-logistic random variable,
then the characteristic function of $y$ is given as
\begin{eqnarray}
\Phi_y(t)&=&\int_{-\infty}^\infty {\rm e}^{it[\mu+\theta
\text{ln}\frac{x^2}{a(q-1)m}]}f(x){\rm d}x\nonumber\\
&=&\frac{\Gamma(\frac{m+1}{2})}{\Gamma(\frac{m}{2})\sqrt{(\pi
m)}}\frac{{\rm
e}^{it\mu}}{[a(q-1)m]^{it\theta}}\int_{-\infty}^\infty
(x^2)^{it\theta}(1+\frac{x^2}{m})^{-(\frac{m+1}{2})}{\rm d}x\nonumber
\end{eqnarray}
Since the integrand is an even function so we can write
\begin{eqnarray}
&=&\frac{\Gamma(\frac{m+1}{2})}{\Gamma(\frac{m}{2})\sqrt{(\pi
m)}}\frac{2~{\rm
e}^{it\mu}}{[a(q-1)m]^{it\theta}}\int_{0}^\infty
(x^2)^{it\theta}(1+\frac{x^2}{m})^{-(\frac{m+1}{2})}{\rm d}x\nonumber\\
&=&\frac{\Gamma(\frac{m+1}{2})}{\Gamma(\frac{m}{2})\sqrt{\pi}}\frac{{\rm e}^{it\mu}}{[a(q-1)m]^{it\theta}}\int_0^\infty
u^{it\theta+\frac{1}{2}-1}(1+u)^{-(\frac{m+1}{2})}{\rm d}u,~~\frac{x^2}{m}=u\nonumber\\
&=&\frac{{\rm
e}^{it\mu}}{[a(q-1)]^{it\theta}}\frac{\Gamma(it\theta+\frac{1}{2})
\Gamma(\frac{m}{2}-it\theta)}{\Gamma(\frac{m}{2})\Gamma(\frac{1}{2})}
\end{eqnarray}
which is the characteristic function of a generalized $q$-logistic
distribution with parameters $(\frac{1}{2},\frac{m+1}{2})$. then by
the uniqueness theorem, the proof is established.
{\section{\bf Skew $q$-Logistic Distribution }}
In this section, we mainly consider two different generalizations of the logistic distribution
by introducing skewness parameters. It may be mentioned that although several skewed
distribution functions exist on the positive real axis, but not many skewed distributions are
available on the whole real line, which are easy to use for data analysis purpose. The main
idea is to introduce the skewness parameter, so that the generalized logistic distribution can
be used to model data exhibiting a unimodal density function having some skewness present
in the data, a feature which is very common in practice.
\vskip.1cm
Recently, skewed distributions have played an important role in the statistical literature
since the pioneering work of Azzalini (1985). He has provided a methodology to introduce
skewness in a normal distribution. Since then a number of papers appeared in this area, see
for example the monograph by Genton (2004) for some recent references.
\vskip.2cm
The first generalization is carried out using the idea of Azzalini (1985) and we name it as
the skew logistic distribution. It is observed that using the same basic principle of Azzalini (1985), the skewness can be easily introduced to the logistic distribution. It has location, scale
and skewness parameters. It is observed that the PDF of the skew logistic distribution can
have different shapes with both positive and negative skewness depending on the skewness
parameter. It has heavier tails than the skew normal distribution proposed by Azzalini (1985).
Although the PDF of the skew logistic distribution is unimodal and log-concave, but the
distribution function, failure rate function and the different moments can not be obtained
in explicit forms. Moreover, it is observed that even when the location and scale parameters
are known, the maximum likelihood estimator of the skewness parameter may not always
exist. Due to this problem, it becomes difficult to use this distribution for data analysis
purposes.
\vskip.2cm
Azzalini (1985) proposed the skew normal distribution, which has the following density function;
\begin{equation}
h(y;\alpha)=2 \Phi (\alpha y) \phi(y);-\infty<y<\infty,
\end{equation}
here $\alpha$ is the skewness parameter, $\phi(\cdot)$ and $\Phi(\cdot)$ are the density function and distribution
function of the standard normal random variable. A motivation of the above model has been
elegantly exhibited by Arnold et al. (1993). Although, Azzalini has been extended the standard
normal distribution function to the form (9), but it has been observed that similar method
can be applied to any symmetric density function. For example, if $f(\cdot)$ is any symmetric
density function defined on $(-\infty, \infty)$ and $F(\cdot)$ is its distribution function, then for any
$\alpha \in(-\infty, \infty)$,
\begin{equation}
2 F (\alpha y) f(y);-\infty<y<\infty,
\end{equation}
is a proper density function and it is skewed if $\alpha\neq0$. This property has been studied
extensively in the literature to study skew-t and skew-Cauchy distributions. Along the same
line we define the skew $q$-logistic distribution with the skewness parameter $\alpha$ as follows. If a
random variable $y$ has the following density function
\begin{equation}\label{eq:1}
f_1(y;\alpha)=2 F (\alpha y) f(y);-\infty<y<\infty,
\end{equation}
then we say that $y$ has a skew$q$-logistic (S$q$L) distribution with skewness parameter $\alpha$. For
brevity we will denote it by $S$q$L(\alpha)$. Clearly (\ref{eq:1}) is a proper density function and $\alpha=0$, corresponds to the standard $q$-logistic distribution.
\begin{center}
\begin{eqnarray*}
\begin{array}{cc}
 \resizebox{6cm}{!}{\includegraphics{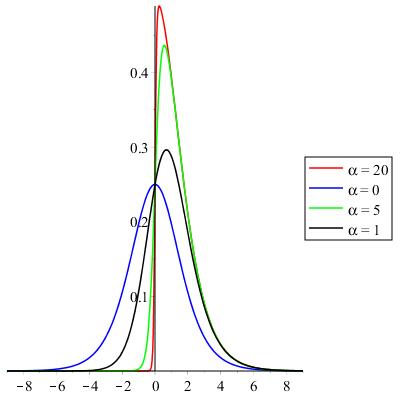}}~~~~~~~&~~~~~~~ \resizebox{6cm}{!}{\includegraphics{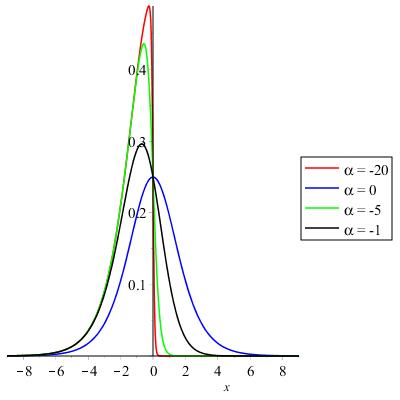}}\\
 \text{Figure 8: {\scriptsize $q=2$
 and $\alpha=0,1,5,20$}}&  \text{Figure 9: {\scriptsize $q=2$
 and $\alpha=0,-1,-5,-20$ }}\\
 \end{array}
\end{eqnarray*}
\end{center}
Above figures, it is clear that $SqL(\alpha)$ is positively skewed when $\alpha$ is positive. It takes
similar shapes on the negative side for $\alpha < 0$. Therefore, $SqL(\alpha)$ can take positive and
negative skewness. As $\alpha$ goes to $\pm \infty$, it converges to the half logistic distribution. Comparing
with the shapes of the skew normal density function of Azzalini (1985), it is clear that $SqL(\alpha)$
produces heavy tailed skewed distribution than the skew normal ones. For large values of
$\alpha$, the tail behaviors of the different members of the $SqL(\alpha)$ family are very similar, which
is apparent from (\ref{eq:1}) also. It is clear from Figure 2 that the tail behaviors of the different
family members of $SqL(\alpha)$ are same for large values of $\alpha$. Some of the properties which are
true for skew normal distribution are also true for skew $q-$logistic distribution.

\end{document}